\documentclass[11pt]{article}
\usepackage{enumerate}
\usepackage{amssymb,a4wide,latexsym,makeidx,epsfig,fleqn}
\usepackage{amsthm}
\usepackage{amsmath}
\usepackage{enumerate}
\newtheorem{theorem}{Theorem}[section]
\newtheorem{remark}[theorem]{Remark}

\newtheorem{lemma}[theorem]{Lemma}

\begin{document}
\textwidth 150mm \textheight 225mm
\title{On the signless Laplacian spectral radius of $C_{4}$-free $k$-cyclic graphs
\thanks{ Supported by the National Natural Science Foundation of China (No. 11171273) and sponsored by the Seed Foundation of Innovation and Creation for Graduate Students in Northwestern Polytechnical University (No. Z2016170)
 \vskip 0.05in}}
\author{{Qi Kong$^1$, Ligong Wang$^2$}\\
{\small Department of Applied Mathematics, School of Science, Northwestern
Polytechnical University,}\\ {\small  Xi'an, Shaanxi 710072,
People's Republic
of China.}\\ {\small $^1$E-mail: kongqixgd@163.com}\\
{\small $^2$E-mail: lgwangmath@163.com}}
\date{}
\maketitle
\begin{center}
\begin{minipage}{120mm}
\vskip 0.3cm
\begin{center}
{\small {\bf Abstract}}
\end{center}
{\small A $k$-cyclic graph is a connected graph of order $n$ and size $n+k-1$. In this paper, we determine the maximal signless Laplacian spectral radius and the corresponding extremal graph among all $C_{4}$-free $k$-cyclic graphs of order $n$. Furthermore, we determine the first three unicyclic, and bicyclic, $C_{4}$-free graphs whose  spectral radius of the signless Laplacian is maximal. Similar results are obtained for the (combinatorial) Laplacian.

\vskip 0.1in \noindent {\bf Key Words}: \ $k$-cyclic graph; $C_{4}$-free; signless Laplacian spectral radius; Laplacian spectral radius. \vskip
0.1in \noindent {\bf AMS Subject Classification (1991)}: \ 05C50, 15A18. }
\end{minipage}
\end{center}

\section{Introduction }
\label{sec:ch6-introduction}

Let $G=(V(G), E(G))$ be a simple graph with vertex set $V(G)$ and edge set $E(G)$. Denote by $v(G)$ the order of $G$ and $e(G)$ the size of $G$, that is to say, $v(G)=|V(G)|$, and $e(G)=|E(G)|$. $\Gamma_{G}(u)=\{v|uv\in E(G)\}$ and $d_{G}(u)=|\Gamma_{G}(u)|$, or simply $\Gamma(u)$ and $d(u)$, respectively. Let $\delta=\delta(G)$ and $\Delta=\Delta(G)$ denote the minimum degree and maximum degree of the graph $G$. Let $X$ and $Y$ be disjoint subsets of $V(G)$. $e(X,Y)$ is the number of edges (in $G$) joining vertices in $X$ to vertices in $Y$ for $G$.  Let $P_{n}$, $S_{n}$, $C_{n}$ and $K_{n}$ be the path, star, cycle and complete graph of order $n$, respectively.

The union of $G_{1}$ and $G_{2}$ is the graph $G_{1}\cup G_{2}$, whose vertex set is $V_{1} \cup V_{2}$ and whose edge set is $E_{1} \cup E_{2}$. $kG$ denotes the union of $k$ copies of $G$. The join of graphs $G_{1}$ and $G_{2}$ is the graph $G_{1}\vee G_{2}$ obtained from $G_{1}\cup G_{2}$ by joining each vertex of $G_{1}$ with every vertex of $G_{2}$. Let $G_{n}^{k}=K_{1}\vee(kK_{2}\cup(n-2k-1)K_{1})$ (see Fig. 1). We use $U_{n}$ to denote the family of all unicyclic graphs of order $n$, and $B_{n}$ to denote the family of all  bicyclic graphs of order $n$.

The matrix $Q(G)=D(G)+A(G)$ is called the signless Laplacian matrix of $G$, where $A(G)$ is the adjacency matrix of $G$ and $D(G)$ is the diagonal matrix of vertex degrees of $G$. The matrix $L(G)=D(G)-A(G)$ is called the Laplacian matrix of $G$. The largest eigenvalue of $A(G)$, $L(G)$ and $Q(G)$ are called the spectral radius, Laplacian spectral radius and signless Laplacian spectral radius (or $Q$-index) of $G$, respectively, and denoted by $\rho(G)$, $\mu(G)$ and $q(G)$, respectively.

The central problem of the classical extremal graph theory is the Tur\'{a}n's Problem:

Problem A. \textit{Given a graph $F$, what is the maximum number of edges of a graph of order $n$, with no subgraph isomorphic to $F$?}

Such problems are well understood nowadays, for example, see the book \cite{Bol} and the survey \cite{BrHo}. Recently, Nikiforov, et al., investigated spectral Tur\'{a}n'
s Problem, namely for the spectral radius $\rho(G)$ of a graph $G$. In this new class of problems the central question is the following one.

Problem B. \textit{Given a graph $F$, what is the maximum spectral radius $\rho(G)$ of a graph $G$ of order $n$, with no subgraph isomorphic to $F$?}

In \cite{BrHo}, when the graph $F$ is the complete graph $K_{r}$, the path or the cycle etc., Nikiforov determines the largest spectral radius of the graph $G$ and their corresponding extremal graphs. The present paper contributes to an even newer trend in extremal graph theory, namely to the study of variations of Promble A for the signless Laplacian spectral radius of graphs, where the central question is the following one.

Problem C. \textit{Given a graph $F$, what is the maximum signless Laplacian spectral radius $q(G)$ of a graph $G$ of order $n$, with no subgraph isomorphic to $F$?}

In \cite{BrSo}, Nikiforov et al. determine the maximum signless Laplacian spectral radius of graphs with  no 4-cycle and 5-cycle. And when the graph $F$ is the complete graph $K_{r}$, the path or the cycle etc., the authors in \cite{Fen, CvSi, FEC, FLZ} determine the largest signless Laplacian spectral radius of the graph $G$ and their corresponding extremal graphs, respectively. He and Guo in \cite{FYI} determine the extremal graph of the signless Laplacian and Laplacian spectral radius among $C_{3}$-free $k$-cyclic graphs of order $n$.

In this paper, we determine the signless Laplacian spectral radius of $C_{4}$-free $k$-cyclic graphs of order $n$ and characterize its extremal graph. Furthermore, we determine the first three signless Laplacian spectral radius of $C_{4}$-free unicyclic graphs of order $n$, and $C_{4}$-free bicyclic graphs of order $n$ with.

\section{Main Lemmas}
In this section, we state some well-known results which will be used in this paper.

\begin{lemma}\label{le:1} (\cite{HoZh})
Let $G$ be a connected graph of order $n$ and $q(G)$ its signless Laplacian spectral radius of $G$. Let $u$, $v$ be two vertices of $G$ and $d(v)$ be the degree of vertex $v$. Assume $v_{1}, v_{2}, \ldots, v_{s}(1\leq s\leq d(v))$ are some vertices of $\Gamma(v)\setminus \Gamma(u)$ and $X=(x_{1}, x_{2}, \ldots, x_{n})^{T}$ is the Perron vector of $Q(G)$, where $x_{i}$ corresponds to the vertex $v_{i}$ $(1\leq i\leq n)$. Let $G^{*}$ be the graph obtained from $G$ by deleting the edges $vv_{i}$ and adding the edges $uv_{i}(1\leq i\leq s)$, If $x_{u}\geq x_{v}$, then $q(G)<q(G^{*})$.
\end{lemma}

\begin{lemma}\label{le:2} (\cite{FeYu,Mer}) For every graph $G$, we have
$$\displaystyle q(G)\leq\max\limits_{u\in V(G)}\{d(u)+\frac{1}{d(u)}\sum\limits_{v\in \Gamma(u)}d(v)\}.$$
If $G$ is connected, equality holds if and only if $G$ is regular or semiregular bipartite.
\end{lemma}

\begin{lemma}\label{le:3} (\cite{TaWa, WXH}) Let G be a simple connected graph on $n$ vertices with maximum degree $\Delta$ and at least one edge. Then

\begin{center}
$\mu(G)\geq\Delta(G)+1$,  $q(G)\geq\Delta(G)+1$,
\end{center}
 where the former equality holds if and only if $\Delta(G)=n-1$, and the latter one holds if and only if $G$ is the star $S_{n}$.
\end{lemma}

\begin{figure}[htb]
\centering
\includegraphics[scale=0.3,width=30mm]{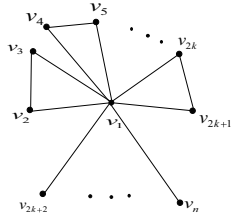}
\caption{The $k$-cyclic graph $G_{n}^{k}$}
\end{figure}
\begin{figure}[htb]
\centering
\includegraphics[scale=0.3,width=100mm]{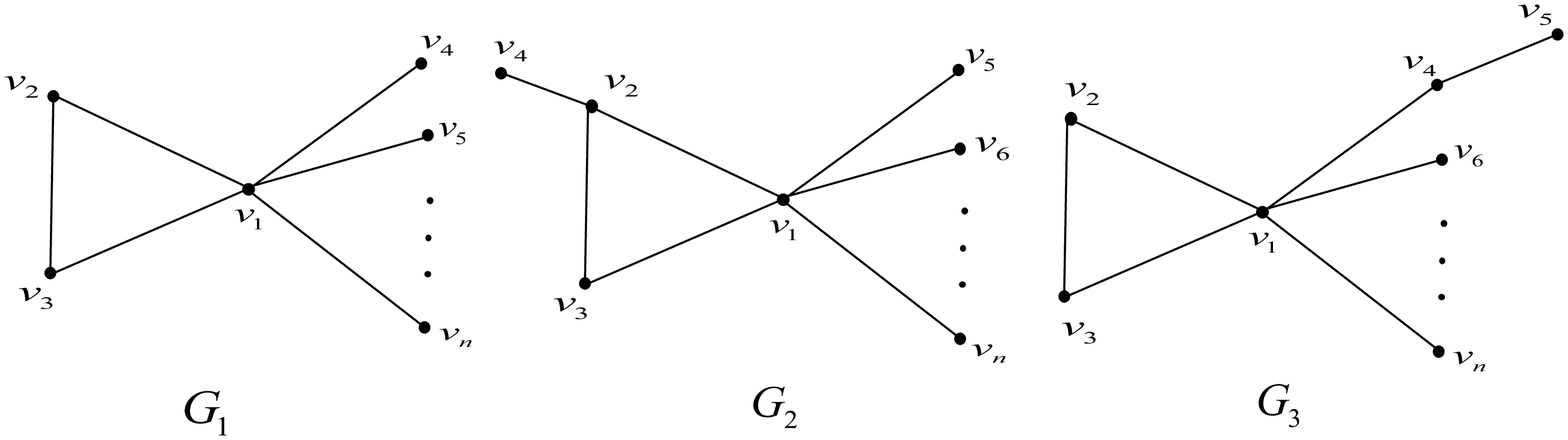}
\caption{The first three unicyclic graphs}
\end{figure}
\begin{figure}[htb]
\centering
\includegraphics[scale=0.3,width=100mm]{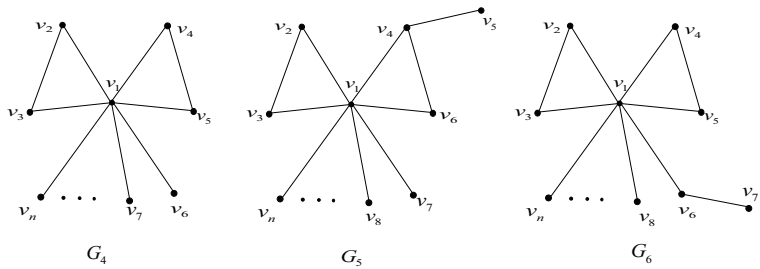}
\caption{The first three bicyclic graphs}
\end{figure}

\begin{lemma}\label{le:4} Let $G_{n}^{k}$, $G_{2}$ and $G_{3}$ be the graphs shown in Figs. 1 and 2. Then $q(G_{n}^{k})$, $q(G_{2})$ and $q(G_{3})$ are the largest roots of the following equations, respectively,
\begin{align*}
&f_{k}(x)\triangleq x^{3}-(n+3)x^{2}+3nx-4k=0,\\
&f_{2}(x)\triangleq x^{5}-(n+5)x^{4}+(6n+3)x^{3}-(9n-1)x^{2}+(3n+8)x-4=0,\\
&f_{3}(x)\triangleq x^{5}-(n+5)x^{4}+(6n+4)x^{3}-(10n-2)x^{2}+(3n+12)x-4=0.
\end{align*}
\end{lemma}

{\bf Proof.}
Let $V(G_{n}^{k})=\{v_{1}, v_{2}, \ldots, v_{n}\}$ and $X=(x_{1}, x_{2}, \ldots, x_{n})^{T}$ be a Perron vector corrosponding to $q(G_{n}^{k})$. By the symmetry of $G_{n}^{k}$, we have

\begin{center}
$x_{2}=x_{3}=\cdots=x_{2k}=x_{2k+1}$,  ~~~$x_{2k+2}=\cdots=x_{n}$.
\end{center}

Let $q=q(G_{n}^{k})$, from the eigenequations for $Q(G_{n}^{k})$ we see that
 \begin{displaymath}
\left\{\
        \begin{array}{ll}
          (q-n+1)x_{1}&=2kx_{2}+(n-2k-1)x_{n},\\
            ~~~~~~(q-3)x_{2}&=x_{1},\\
            ~~~~~~(q-1)x_{n}&=x_{1}.\\
        \end{array}
      \right.
\end{displaymath}
Since $X=(x_{1}, x_{2}, \ldots, x_{n})^{T}$ is an eigenvector corresponding to $q(G_{n}^{k})$, so $X\neq0$. Then
\begin{displaymath}
\left|\begin{array}{ccc}
 \ q-n+1  & -2k & -n+2k+1  \\
     -1  &  q-3 & 0 \\
     -1 &  0  & q-1\\
\end{array}\right|=0.
\end{displaymath}
So $q$ is the largest root of the following equation
\begin{displaymath}
\left|\begin{array}{ccc}
 \ x-n+1  & -2k & -n+2k+1  \\
     -1  &  x-3 & 0 \\
     -1 &  0  & x-1\\
\end{array}\right|=0.
\end{displaymath}
Consequently, $q$ is the largest root of the following equation

\begin{center}
$f_{k}(x)\triangleq x^{3}-(n+3)x^{2}+3nx-4k=0$.
\end{center}
Using the same method, we obtain $q(G_{2})$ and $q(G_{3})$ are the largest roots of the following equations, respectively,
\begin{align*}
&f_{2}(x)\triangleq x^{5}-(n+5)x^{4}+(6n+3)x^{3}-(9n-1)x^{2}+(3n+8)x-4=0.\\
&f_{3}(x)\triangleq x^{5}-(n+5)x^{4}+(6n+4)x^{3}-(10n-2)x^{2}+(3n+12)x-4=0.
\end{align*}
This completes the proof. $\square$

\section{Main results}
Let $G$ be a $C_{4}$-free $k$-cyclic graph of order $n$. If $k=0$ then $G$ is a tree. In \cite{Zhan, Ord}, the authors determined the first eight Laplacian spectral radius of trees of order $n$. For a bipartite graph $G$, by \cite{WXH}, we know that  $L(G)$ and $Q(G)$ have the same eigenvalues. A tree is a bipartite graph, so the results that are obtained by \cite{Zhan, Ord} hold also for the signless Laplacian spectral radius of trees of order $n$.  Therefore, in what follows, assume that $k\geq1$.

\begin{theorem}\label{th:1} Let $k\geq 1, n\geq 2k+2$ and let $G$ be a $C_{4}$-free $k$-cyclic graph of order $n$. Then
$$q(G)\leq q(G_{n}^{k}),$$
with equality if and only if $G=G_{n}^{k}$, where $q(G_{n}^{k})$ is the largest root of the equation
$$x^{3}-(n+3)x^{2}+3nx-4k=0.$$
\end{theorem}

{\bf Proof.}
 By Lemma 2.3, we have

\begin{center}
$q(G_{n}^{k})> \Delta(G_{n}^{k})+1=n-1+1=n$.
\end{center}
Since $n\geq 2k+2$ and $G$ is a $C_{4}$-free $k$-cyclic graph of order $n$, we have $\Delta(G)\leq n-1$. If $\Delta(G)=n-1$, it is easy to see that $G$ must be $G_{n}^{k}$. If $G\neq G_{n}^{k}$, then $\Delta(G)\leq n-2$.
By Lemma 2.2, let $w$ be a vertex of $G$ such that

\begin{center}
$\displaystyle d(w)+\frac{1}{d(w)}\sum\limits_{i\in \Gamma(w)}d(i)=\max\limits_{u\in V(G)}\{d(u)+\frac{1}{d(u)}\sum\limits_{v\in \Gamma(u)}d(v)\}.$
\end{center}
Then

\begin{center}
$\displaystyle q(G)\leq d(w)+\frac{1}{d(w)}\sum\limits_{i\in \Gamma(w)}d(i),$
\end{center}
and $1\leq d(w)\leq\Delta(G)\leq n-2$.

If $d(w)=1$ we obtain
$$\displaystyle q(G)\leq d(w)+\frac{1}{d(w)}\sum\limits_{i\in \Gamma(w)}d(i)\leq 1+\frac{\Delta}{1}\leq n-1<q(G_{n}^{k}).$$
If $2\leq d(w)\leq n-2$, since $G$ is $C_{4}$-free, every vertex $v\in V(G)\setminus\Gamma(w)$  has at most one neighbor in $\Gamma(w)$.
Then we have
$$e(\Gamma(w),V(G)\setminus\Gamma(w))\leq d(w)+| V(G)\setminus(\Gamma(w)\cup\{w\})|=d(w)+n-d(w)-1=n-1.$$
Thus
$$\displaystyle \sum\limits_{i\in \Gamma(w)}d(i)\leq d(w)+e(\Gamma(w),V(G)\setminus\Gamma(w))\leq d(w)+n-1.$$

Moreover, $$\displaystyle q(G)\leq d(w)+\frac{1}{d(w)}\sum\limits_{i\in \Gamma(w)}d(i)\leq 1+d(w)+\frac{n-1}{d(w)}.$$
Since the function
$$\displaystyle f(x)=x+\frac{n-1}{x}$$
is convex for $x>0$, its maximum in any closed interval is attained at one of the ends of this interval.
If $2\leq d(w)\leq n-2$, then

\begin{center}
$\displaystyle q(G)\leq d(w)+\frac{1}{d(w)}\sum\limits_{i\in \Gamma(w)}d(i)\leq 1+\max\{2+\frac{n-1}{2},n-2+\frac{n-1}{n-2}\}\leq n\leq q(G_{n}^{k})$.
\end{center}
From the above all, we obtain $q(G)\leq q(G_{n}^{k})$, with equality if and only if $G=G_{n}^{k}$. Then from Lemma 2.4,  we know that $q(G_{n}^{k})$ is the largest root of the polynomial
$$x^{3}-(n+3)x^{2}+3nx-4k=0.$$
This completes the proof.  $\square$

\begin{theorem}\label{th:2} Let $U_{n}$ be the set unicyclic graphs of order $n$, and $n\geq6$, $G_{i}\in U_{n}$ for $i=1,2,3$. Then for any $G\in U_{n}$ and $G\neq G_{i}  (i=1,2,3 )$, we have
$$q(G)<q(G_{3})<q(G_{2})<q(G_{1}),$$
where $G_{1}, G_{2}$ and $G_{3}$ are the graphs shown in Fig. 2.
\end{theorem}

{\bf Proof.} For any $G\in U_{n}$, from Theorem 3.1, if $G\neq G_{1}$, then $q(G)<q(G_{1})$. Especially,  $q(G_{j})<q(G_{1})$ for $j=2,3$. Obviously, $G_{2}$ and $G_{3} $ are all the unicyclic graphs with $\Delta=n-2$ in $U_{n}$. Then for any $G\in U_{n}$, if $G\neq G_{i} (i=1,2,3)$, then $\Delta(G)\leq n-3$. By Lemma 2.3, $q(G_{j})> \Delta(G_{j})+1=n-2+1=n-1$  for $j=2,3$. From Lemma 2.2,
$$\displaystyle q(G)\leq\max\limits_{u\in V(G)}\{d(u)+\frac{1}{d(u)}\sum\limits_{v\in \Gamma(u)}d(v)\}.$$
Similar to the proof of Theorem 3.1, we can get $q(G)\leq n-1<q(G_{i})$.

From Lemma 2.4,  when $x\geq q(G_{3})>n-1$ and $n\geq6$, we have
 \begin{align*}
F(x)=f_{3}(x)-f_{2}(x)&=x^{3}-(n-1)x^{2}+4x\\
                      &=(x-n+1)^{3}+2(n-1)(x-n+1)^{2}\\
                      &+(n^{2}-2n+5)(x-n+1)+4(n-1)>0
\end{align*}
So, $f_{2}(q(G_{3}))<0$, then $q(G_{2})>q(G_{3})$.

From the above all, for any $G\in U_{n}$ and $G\neq G_{i},(i=1,2,3 )$, we have
$$
q(G)<q(G_{3})<q(G_{2})<q(G_{1}).
$$
This completes the proof. $\square$

\begin{theorem}\label{th:3} Let $B_{n}$ be the set bicyclic graphs of order $n$, and $n\geq 8$, $G_{i}\in B_{n}$ for $i=4, 5, 6$. Then for any $G\in B_{n}$ and $G\neq G_{i} (i=4,5,6 )$, we have
$$q(G)<q(G_{6})<q(G_{5})<q(G_{4}),$$
where $G_{4}, G_{5}$ and $G_{6}$ are shown in Fig. 3.
\end{theorem}

{\bf Proof.} For any $G\in B_{n}$, from Theorem 3.1, if $G\neq G_{4}$, then $q(G)<q(G_{4})$. Especially, $q(G_{j})<q(G_{4}) $for $j=5,6$. Obviously, $G_{5}$ and $G_{6}$ are all the bicyclic graphs with $\Delta=n-2$ in $B_{n}$, then for any $G\in B_{n}$, if $G\neq G_{i}$ $(i=4,5,6)$, then $\Delta(G)\leq n-3$. By Lemma 2.3, $q(G_{j})> \Delta(G_{j})+1=n-2+1=n-1$ for $j=5,6$.

From Lemma 2.2,
$$\displaystyle q(G)\leq\max\limits_{u\in V(G)}\{d(u)+\frac{1}{d(u)}\sum\limits_{v\in \Gamma(u)}d(v)\}.$$
Similar to the proof of Theorem 3.1, we can get $q(G)\leq n-1<q(G_{i})$.

Note that
$$G_{5}=G_{6}-v_{4}v_{5}+v_{5}v_{6}=G_{6}-v_{6}v_{7}+v_{4}v_{7}.$$
Applying Lemma 2.1 for the vertices $v_{4}$ and $v_{6}$ of the graph $G_{6}$, we have
$$q(G_{6})<q(G_{5}).$$
From the above all, for any $G\in B_{n}$ and $G\neq G_{i},(i=4,5,6 )$,
we have
$$q(G)<q(G_{6})<q(G_{5})<q(G_{4}).$$
This completes the proof. $\square$

\begin{remark}
Lemma 2.3 is also true for the Laplacian special radius of graphs. From the proof of Theorem 3.1, we know that Theorem 3.1 also holds for Laplacian spectral radius of graphs. From \cite{WXH}, we know $\mu(G)\leq q(G)$, the equality holds if and only if $G$ is a bipartite graph. And from the proofs of Theorems 3.2 and 3.3, we know that Theorems 3.2 and 3.3 also hold for Laplacian special radius of graphs.
\end{remark}

\end{document}